\documentclass{jctart08}

\usepackage{cite}
\usepackage{amsmath}
\usepackage{amssymb}
\usepackage{amsthm}
\theoremstyle{plain}
\newtheorem{Theorem}{Theorem}[section] %
\newtheorem{Lemma}{Lemma}[section]

\theoremstyle{definition}
\newtheorem{Def}{Definition}[section]

\renewcommand{\thetable}{\thesection.\arabic{table}}

\addto\captionsrussian{\def\refname{References}}
\newenvironment{Proof} 
{\par\noindent{\it Proof of}} 
{\hfill$\vspace{5mm}\scriptstyle\blacksquare$} 

\numberwithin{equation}{section} 
\numberwithin{figure}{section} 
\numberwithin{table}{section} 

\begin{document}

\setcounter{page}{1}

\markboth{M.I. Isaev}{Exponential instability in the Gel'fand inverse problem on the energy intervals}

\title{Exponential instability in the Gel'fand inverse problem on the energy intervals}
\date{}
\author{ 
\textbf{M.I. Isaev}\\
Moscow Institute of Physics and Technology,
141700 Dolgoprudny, Russia\\
Centre de Math\'ematiques Appliqu\'ees, Ecole Polytechnique,
91128 Palaiseau, France\\\
e-mail: \tt{isaev.m.i@gmail.com}}

\maketitle

{\bf Abstract}
\begin{abstract}
		We consider the Gel'fand inverse problem and continue studies of [Mandache,2001]. We show that the Mandache-type instability 
		remains valid even in the case of 
		Dirichlet-to-Neumann 
		map given on the energy intervals. These instability results show, in particular, that
		the logarithmic stability estimates of [Alessandrini,1988], [Novikov, Santacesaria,2010] and 
		especially of [Novikov,2010] are optimal (up to the value of the exponent).  
\end{abstract}

\section{Introdution}
We consider the  Schr\"odinger equation
\begin{equation}\label{eq} 
	-\Delta \psi  + v(x)\psi = E\psi, \ \  x \in  D,
\end{equation}
where
\begin{equation}\label{eq_c}
 D \text{ is an open bounded domain in } \mathbb{R}^d,\ d \geq 2,\  \partial D \in C^2,\ v \in L^{\infty}(D).
\end{equation}
Consider the map $\Phi(E)$ such that
\begin{equation}\label{def_phi}
	\Phi(E) (\psi|_{\partial D}) = \frac{\partial \psi}{\partial \nu}|_{\partial D}.
\end{equation}
for all sufficiently regular solutions $\psi$ of (\ref{eq}) in $\bar{D} = D \cup \partial D,$ where $\nu$ is the outward normal to $\partial D$. Here we assume also that
\begin{equation}\label{condition}
	\text{$E$ is not a Dirichlet eigenvalue for  operator $-\Delta + v$ in $D$.} 
\end{equation}
The map $\Phi(E)$ is called the Dirichlet-to-Neumann map and is considered
as boundary measurements.

We consider the following inverse boundary value problem for equation (\ref{eq}). 

{\bf Problem 1.1.}
 Given $\Phi(E)$ on the union of the energy intervals $S = \bigcup\limits_{j=1}\limits^{K} I_j$, find $v$. 
 
 Here we suppose that condition (\ref{condition}) is fulfilled for any $E \in S$.

This problem can be considered as the Gel’fand inverse boundary value problem for the Schr\"odinger equation on the energy intervals 
(see \cite{Gelfand1954}, \cite{Novikov1988}). 

Problem 1.1 includes, in particular, the following questions: (a) uniqueness, (b) reconstruction,
(c) stability.

Global uniqueness for Problem 1.1 was obtained for the first time by Novikov (see Theorem 5.3 in \cite{Henkin1987}).  Some global reconstruction 
method for Problem 1.1  was proposed for the first time in 
\cite{Henkin1987} also. Global uniqueness theorems and global reconstruction methods in the case of fixed energy
were given for the first time in \cite{Novikov1988} in dimension $d \geq 3$ and in \cite{Bukhgeim2008} in dimension $d=2$.

Global stability estimates for Problem 1.1 were given for the first time in \cite{Alessandrini1988} 
in dimension $d \geq 3$ and in \cite{NS} in dimension $d=2$. The Alessandrini result
of \cite{Alessandrini1988} was recently improved by Novikov in \cite{Novikov2010}.
In the case of fixed energy, Mandache showed in \cite{Mandache2001} that these logarithmic 
stability results are optimal (up to the value of the exponent). Mandache-type instability estimates for inverse
inclusion and scattering problems are given in \cite{Cristo2003}.

	 In the present work we extend studies of Mandache to the case of Dirichlet-to-Neumann map given on the energy intervals. The stability estimates and our instability results for Problem 1.1 are presented and discussed in Section 2. In Section 5 we prove the main results, using a ball packing and covering by ball arguments. In Section 3 we prove some
 basic properties of the Dirichlet-to-Neumann map, using some Lemmas about the Bessel functions wich we proved in Section 6.
	 
\section{Stability estimates and main results}

As in \cite{Novikov2010} we assume for simplicity that
\begin{equation}\label{assumption}
	\begin{array}{l}
	 D \text{ is an open bounded domain in } \mathbb{R}^d,\  \partial D \in C^2,\\
	 v \in W^{m,1}(\mathbb{R}^d) \text{ for some } m > d, \ \mbox{supp}\, v \subset D,\  d \geq 2,
	\end{array}
\end{equation}
where
\begin{equation}
	W^{m,1}(\mathbb{R}^d) = \{v:\  \partial^J v \in L^1(\mathbb{R}^d),\  |J| \leq m \},\ m\in \mathbb{N}\cup 0,
\end{equation}
where 
\begin{equation}
J \in (\mathbb{N}\cup 0)^d,\ |J| = \sum\limits_{i=1}\limits^{d}J_i,\ \partial^J v(x) 
= \frac{\partial^{|J|} v(x)}{\partial x_1^{J_1}\ldots \partial x_d^{J_d}}.
\end{equation}
Let
\begin{equation}
	||v||_{m,1} = \max\limits_{|J|\leq m} ||\partial^J v||_{L^1(\mathbb{R}^d)}.
\end{equation}
We recall that if $v_1$, $v_2$ are potentials satisfying (\ref{condition}),(\ref{def_phi}), where $E$ and $D$ are fixed, then
\begin{equation}
	\Phi_1 - \Phi_2 \text{ is a compact operator in } L^{\infty}(\partial D),
\end{equation}
where $\Phi_1$, $\Phi_2$ are the DtN maps for $v_1$, $v_2$ respectively, see \cite{Novikov1988}. Note also that 
(\ref{assumption}) $\Rightarrow$ (\ref{eq_c}).

\begin{Theorem}[variation of the result of \cite{Alessandrini1988}, see \cite{Novikov2010}]\label{TheoremA}
Let conditions (\ref{condition}), (\ref{assumption}) hold for potentials $v_1$ and $v_2$, where $E$ and $D$ are fixed,
$d \geq 3$. Let $||v_j||_{m,1} \leq N,\  j = 1,2$, for some $N > 0$. Let $\Phi_1,\ \Phi_2$ denote DtN maps for
$v_1$, $v_2$ respectively. Then
\begin{equation}\label{estimation_A}
	||v_1 - v_2||_{L^\infty(D)} \leq 
	c_1 (\ln(3+||\Phi_1 - \Phi_2||^{-1}))^{-\alpha_1},
\end{equation}
where $c_1 = c_1(N,D,m),\ \alpha_1 = (m-d)/m,\  ||\Phi_1 - \Phi_2|| = 
||\Phi_1 - \Phi_2||_{L^\infty(\partial D)\rightarrow L^\infty(\partial D)}$.
\end{Theorem}
An analog of stability estimate of \cite{Alessandrini1988} for $d = 2$ is given in \cite{NS}.

	A disadvantage of estimate (\ref{estimation_A}) is that
\begin{equation}\label{alpha_A}
	\alpha_1 < 1 \text{ for any } m > d \text{ even if $m$ is very great.} 
\end{equation} 

\begin{Theorem}[the result of \cite{Novikov2010}]\label{TheoremN}
Let the assumptions of Theorem \ref{TheoremA} hold. Then
\begin{equation}\label{estimation_N}
	||v_1 - v_2||_{L^\infty(D)} \leq 
	c_2 (\ln(3+||\Phi_1 - \Phi_2||^{-1}))^{-\alpha_2},
\end{equation}
where $c_2 = c_2(N,D,m),\ \alpha_2 = m-d,\  ||\Phi_1 - \Phi_2|| = 
||\Phi_1 - \Phi_2||_{L^\infty(\partial D)\rightarrow L^\infty(\partial D)}$.
\end{Theorem}

	A principal advantage of estimate (\ref{estimation_N}) in comparison with (\ref{estimation_A}) is that
\begin{equation}\label{alpha_N}
	\alpha_2 \rightarrow +\infty \text{ as } m \rightarrow + \infty,
\end{equation} 
in contrast with (\ref{alpha_A}). 
Note that strictly speaking Theorem \ref{TheoremN} was proved in \cite{Novikov2010}
for $E = 0$ with the condition that $\mbox{supp}\, v \subset D,$ so we cant make use of substitution $v_E = v-E$, since  
condition $\mbox{supp}\, v_E \subset D$  does not hold.

We would like to mention that, under the assumptions of Theorems \ref{TheoremA} and \ref{TheoremN},
according to the Mandache results of \cite{Mandache2001},  estimate (\ref{estimation_N}) can not hold with
$\alpha_2 > m(2d - 1)/d$ for real-valued potentials and with $\alpha_2 > m$ for complex potentials.

	As in \cite{Mandache2001} in what follows we fix $D = B(0,1)$, where $B(x,r)$
is the open ball of radius $r$ centred at $x$.  We fix an orthonormal basis in $L^2(S^{d-1}) = L^2(\partial D)$
\begin{equation}
	\begin{array}{l}
	\displaystyle
		\{f_{jp} : j \geq 0;\  1 \leq p \leq p_j \}, \\
		\text{$f_{jp}$ is a spherical harmonic of degree $j$,}
	\end{array}
\end{equation}
   where $p_j$  is the dimension of the space of spherical harmonics of order $j$,
\begin{equation}
	p_j = \binom {j+d-1} {d-1} - \binom {j+d-3} {d-1},
\end{equation}
	where 
\begin{equation}
	\binom {n} {k} = 
		\frac{n(n-1)\cdots(n-k+1)}{k!} \ \ \ \text{ for $n \geq 0$} 
\end{equation}	
and
\begin{equation}
	\binom {n} {k} = 
		0 \ \ \ \text{ for $n < 0$.} 
\end{equation}
The precise choice of $f_{jp}$ is irrelevant for our purposes. Besides orthonormality, we only need
$f_{jp}$ to be the restriction of a homogeneous harmonic polynomial of degree $j$ to the sphere
and so $|x|^j f_{jp}(x/|x|)$ is harmonic. In the Sobolev spaces $H^s(S^{d-1})$ we will use the norm
\begin{equation}
	|| \sum_{j,p} c_{jp}f_{jp}||^2_{H^s} = \sum_{j,p}(1+j)^{2s}|c_{jp}|^2.
\end{equation} 
The notation $(a_{jpiq} )$ stands for a multiple sequence.
We will drop the subscript 
\begin{equation}
	0 \leq j,\  1 \leq p \leq p_j ,\  0 \leq i,\  1 \leq q \leq p_i.  
\end{equation}
We use notations: $|A|$ is the cardinality of a set $A$, $[a]$ is the integer part of real number $a$
and 
$(r, \omega) \in \mathbb{R}_+ \times S^{d-1}$ are polar coordinates for $r\omega=x\in\mathbb{R}^d$.

	The interval $I = [a,b]$ will be referred as $\sigma$-regular interval if for any potential 
$v \in L^{\infty}(D)$	with $||v||_{L^{\infty}(D)}\leq \sigma$ and any $E\in I$ condition (\ref{condition}) is fulfilled. 
Note that for any $E\in I$ and any Dirichlet eigenvalue $\lambda$ for operator $-\Delta$ in $D$ we have that
\begin{equation}\label{lambdaE}
	|E-\lambda|\geq\sigma.
\end{equation}			
It follows from the definition of $\sigma$-regular interval, taking $v\equiv E-\lambda$.
\begin{Theorem}\label{main}
	For $\sigma > 0$ and dimension $d \geq 2$ consider the union $S = \bigcup\limits_{j=1}\limits^{K} I_j$ of $\sigma$-regular intervals. 
	Then for any $m > 0$ and any $s \geq 0$ there is a constant $\beta > 0$, such that
	for any $\epsilon \in (0,\sigma/3)$ and $v_0\in C^m(D)$ with $||v_0||_{L^{\infty}(D)}\leq \sigma/3$ and $\mbox{supp}\, v_0 \subset B(0,1/3)$ 
	there are real-valued potentials $v_1, v_2 \in C^{m}(D)$, also supported in $B(0,1/3)$, such that
	\begin{equation}\label{maineq}
	\begin{array}{l}
	\displaystyle
	 \sup\limits_{E\in S}\Big(||\Phi_1(E) - \Phi_2(E)||_{H^{-s} \rightarrow H^{s}}\Big) \leq \exp\Big(-\epsilon^{-\frac{1}{2m}} \Big), \\
	 ||v_1-v_2||_{L^\infty(D)} \geq \epsilon, \\
		||v_i-v_0||_{C^m(D)} \leq \beta, \ \  i=1,2, \\
		||v_i-v_0||_{L^\infty(D)} \leq \epsilon, \ \ \  i=1,2, \\
	\end{array}
\end{equation}
where $\Phi_1(E)$, $\Phi_2(E)$ are the DtN maps for $v_1$ and $v_2$ respectively.
\end{Theorem}
\noindent
{\bf Remark 2.1.}
 	We can allow $\beta$ to be arbitrarily small in Theorem \ref{main}, if we 
 	require $\epsilon \leq \epsilon_0$ and replace the right-hand side in the instability estimate
 	by $\exp(-c\epsilon^{-\frac{1}{2m}} )$, with $\epsilon_0 > 0$ 	and $c > 0$, depending on $\beta.$

In addition to Theorem \ref{main}, 
we consider explicit instability example with a complex potential given by Mandache in \cite{Mandache2001}. 
We show that it gives exponential instability even in case of Dirichlet-to-Neumann map given on the energy intervals.
Consider the cylindrical variables $(r_1, \theta, x') \in \mathbb{R_+}\times\mathbb{R}/2\pi\mathbb{Z}\times \mathbb{R}^{d-2}$,
with $x' = (x_3, \ldots, x_d)$, $r_1\cos\theta = x_1$ and $r_1\sin\theta = x_2$. Take $\phi \in C^{\infty}(\mathbb{R}^2)$
with support in $B(0,1/3)\cap \{x_1 > 1/4\}$ and with $||\phi||_{L^\infty}=1$. 

\begin{Theorem}\label{complex}
	For $\sigma > 0$, $m > 0$, integer $n > 0$ and dimension $d \geq 2$ consider the union $S = \bigcup\limits_{j=1}\limits^{K} I_j$ of $\sigma$-regular intervals
	and  define the complex potential
	\begin{equation}
		v_{nm}(x) = \frac{\sigma}{3} n^{-m} e^{in\theta}\phi(r_1,|x'|). 
	\end{equation}
	Then $||v_{mn}||_{L^\infty(D)} = \frac{\sigma}{3} n^{-m} $ and for every $s \geq 0$ and $m > 0$ there are constants $c$, $c'$ 
	such that $||v_{mn}||_{C^m(D)} \leq c$ and for every $n$
	\begin{equation} 
	\sup\limits_{E\in S}\Big(||\Phi_{mn}(E) - \Phi_0(E)||_{H^{-s} \rightarrow H^{s}}\Big) \leq c'2^{-n/4}, 
	\end{equation}
where $\Phi_{mn}(E)$, $\Phi_0(E)$ are the DtN maps for $v_{mn}$ and $v_0\equiv0$ respectively.	
\end{Theorem}

\noindent 
In some important sense, this is stronger than Theorem \ref{main}.  Indeed, if we take $\epsilon = \frac{\sigma}{3} n^{-m}$ we obtain 
(\ref{maineq}) with $\exp(-C\epsilon^{-1/m})$ in the right-hand side. An explicit real-valued counterexample should be 
difficult to find. This is due to nonlinearity of the map $v \rightarrow \Phi$.

\noindent 
{\bf Remark 2.2.}
Note that for sufficient large $s$ one can see that 
 	\begin{equation}
 		||\Phi_1-\Phi_2||_{L^{\infty}(\partial D) \rightarrow L^{\infty}(\partial D)} 
 		\leq C||\Phi_1-\Phi_2||_{H^{-s} \rightarrow H^{s}}.
 	\end{equation}
 	So Theorem \ref{main} and Theorem \ref{complex} imply, in particular, that the estimate
 	\begin{equation}
 	||v_1 - v_2||_{L^\infty(D)} \leq 
	c_3 \sup_{E \in S} 
	\left(\ln(3+||\Phi_1(E) - \Phi_2(E)||^{-1})\right)^{-\alpha_3},
	\end{equation}
	where $c_3 = c_3(N,D,m,S)$ and $||\Phi_1(E) - \Phi_2(E)|| = 
||\Phi_1(E) - \Phi_2(E)||_{L^\infty(\partial D)\rightarrow L^\infty(\partial D)}$, can not hold
 	with $\alpha_3 > 2m$ for real-valued potentials and with $\alpha_3 > m$ for complex potentials.
 	Thus Theorem \ref{main} and Theorem \ref{complex} show optimality of
 	logarithmic stability results of Alessandrini and Novikov in considerably stronger sense that results of Mandache.

\section{Some basic properties of Dirichlet-to-Neumann map}

 We continue to consider $D = B(0,1)$ and also to use polar coordinates 
$(r, \omega) \in \mathbb{R}_+ \times S^{d-1}$, with $x = r\omega$.
 Solutions of equation $-\Delta \psi = E\psi$ in $D$ can be expressed by the Bessel functions $J_\alpha$ and $Y_\alpha$ with integer or half-integer order $\alpha$, see definitions of Section 6.  Here we state some Lemmas about these functions 
 (Lemma \ref{Lemma_solution}, Lemma \ref {Lemma_ring} and Lemma \ref{Bessel}). 

\begin{Lemma}\label{Lemma_solution}
	Suppose $k \neq 0$ and $k^2$ is not a Dirichlet eigenvalue for  operator $-\Delta$ in $D$. Then
	\begin{equation}\label{psi_0}
		\psi_0(r,\omega) = {r^{-\frac{d-2}{2}}}\frac{J_{j+\frac{d-2}{2}}(kr)}{J_{j+\frac{d-2}{2}}(k)}f_{jp}(\omega) 
	\end{equation}	
	is the solution of equation (\ref{eq}) with $v \equiv 0$, $E = k^2$ and boundary condition $\psi|_{\partial D} = f_{jp}.$
\end{Lemma}
\noindent
{\bf Remark 3.1.}
	Note that the assumptions of Lemma \ref{Lemma_solution} imply  $J_{j+\frac{d-2}{2}}(k) \neq 0$. 

\begin{Lemma}\label{Lemma_ring}
	Let the assumptions of Lemma \ref{Lemma_solution} hold. Then system of functions
		\begin{equation}
			\left\{\psi_{jp}(r,\omega) = R_j(k,r)f_{jp}(\omega) : j \geq 0; 1\leq p\leq p_j \right\}, 
		\end{equation}
		where 
		\begin{equation}
			R_j(k,r) = {r^{-\frac{d-2}{2}}} \Big(Y_{j+\frac{d-2}{2}}(kr)J_{j+\frac{d-2}{2}}(k) - 
			J_{j+\frac{d-2}{2}}(kr)Y_{j+\frac{d-2}{2}}(k) \Big),
		\end{equation}
	is complete orthogonal system (in the sense of $L_2$) in the space of solutions of equation (\ref{eq}) in 
	$D' = B(0,1)\setminus B(0,1/3)$ with $v \equiv 0$, $E = k^2$ and boundary condition $\psi|_{r=1} = 0.$ 
\end{Lemma}

\begin{Lemma}\label{Bessel}
	For any $C>0$ and integer $d \geq 2$ there is a constant $N>3$ depending on $C$ such that for any integer $n\geq N$ and any $|z|\leq C$
		\begin{equation}\label{Jn}
			\frac{1}{2}\frac{(|z|/2)^\alpha}{\Gamma(\alpha+1)} \leq |J_\alpha(z)| \leq \frac{3}{2}\frac{(|z|/2)^\alpha}{\Gamma(\alpha+1)},
		\end{equation}
		\begin{equation}\label{Jn'}
			|J'_\alpha(z)| \leq {3}\frac{(|z|/2)^{\alpha-1}}{\Gamma(\alpha)},
		\end{equation}
		\begin{equation}\label{Yn}
			\frac{1}{2\pi}(|z|/2)^{-\alpha}\Gamma(\alpha) \leq |Y_\alpha(z)| \leq \frac{3}{2\pi}(|z|/2)^{-\alpha}\Gamma(\alpha)
		\end{equation}
		\begin{equation}\label{Yn'}
			|Y'_\alpha(z)| \leq \frac{3}{\pi}(|z|/2)^{-\alpha-1}\Gamma(\alpha+1)
		\end{equation}
		where $'$ denotes derivation with respect to $z$, $\alpha = n+\frac{d-2}{2}$ and $\Gamma(x)$ is the Gamma function.
\end{Lemma}  

Proofs of Lemma \ref{Lemma_solution}, Lemma \ref {Lemma_ring} and Lemma \ref{Bessel} are given in Section 6.

\begin{Lemma}\label{Lemma1}
	Consider a compact $W \subset \mathbb{C}$.  
Suppose, that   $v$ is bounded, $\mbox{supp}\, v \subset B(0,1/3)$ and condition (\ref{condition}) is fulfilled for any $E \in W$ and 
potentials $v$ and $v_0$, where $v_0 \equiv 0.$ Denote $\Lambda_{v,E} = \Phi(E) - \Phi_0(E).$ 
Then there is a constant $\rho = \rho(W,d)$, such that for any 
$0\leq j, 1\leq p\leq p_j$, $0\leq i, 1\leq q\leq p_i$, we have
\begin{equation}\label{eqlemma1}
\left|\left\langle \Lambda_{v,E} f_{jp},f_{iq}\right\rangle\right|
\leq
\rho \,2^{-\max(j,i)}||v||_{L^\infty(D)}||(-\Delta+v-E)^{-1}||_{L^2(D)},
\end{equation}
where  $\Phi(E)$, $\Phi_0(E)$ are the DtN maps for $v$ and $v_0$ respectively and 
$(-\Delta+v-E)^{-1}$ is considered with the Dirichlet boundary condition.
\end{Lemma}
\begin{Proof} {\it Lemma \ref{Lemma1}. }
	For simplicity we give first a proof under the additional assumtions that $0 \notin W$ and there is a holomorphic germ $\sqrt{E}$ for $E\in W$.  
	Since $W$ is compact there is $C>0$ such that for any $z\in W$ we have $|z| \leq C$. We take $N$ from Lemma \ref{Bessel} for this $C$.  
	We fix indeces $j, p$. Consider solutions  $\psi(E)$, $\psi_0(E)$ of equation (\ref{eq})   with $E \in W$, 
	boundary condition $\psi|_{\partial D} = f_{jp}$ and potentials $v$ and $v_0$ respectively. 
	Then $\psi(E) - \psi_0(E)$ has zero boundary values, so it is domain of $-\Delta + v - E$, and since 
	\begin{equation}
	(-\Delta + v - E)\left(\psi(E) - \psi_0(E)\right) = -v\psi_0(E) \  \text{ in $D$}, 
	\end{equation}
	we obtain that
	\begin{equation}\label{analytic_eq}
		\psi(E) - \psi_0(E) = -(-\Delta + v - E)^{-1}v\psi_0(E).
	\end{equation}
	If $j \geq N$ from Lemma \ref{Lemma_solution} and Lemma \ref{Bessel} we have that
	\begin{equation}
		\begin{array}{l}
		\displaystyle
		||\psi_0(E)||_{L^2(B(0,1/3))}^2= 
		||f_{jp}||_{L^2(S^{d-1})}^2
		\int_0^{1/3}
		\left|r^{-\frac{d-2}{2}}
		\frac
		{J_{j+\frac{d-2}{2}}(\sqrt{E}\,r)} 
		{J_{j+\frac{d-2}{2}}(\sqrt{E})}\right|^2  r^{d-1}dr\leq
		\\\displaystyle
		\leq \int_0^{1/3} 
		\Bigg(\frac{3}{2}\frac{(|E|^{1/2}r/2)^{j+\frac{d-2}{2}}}{\Gamma(j+\frac{d-2}{2}+1)} \Bigg)^2\Big/ 
		\Bigg( \frac{1}{2}\frac{(|E|^{1/2}/2)^{j+\frac{d-2}{2}}}{\Gamma(j+\frac{d-2}{2}+1)} \Bigg)^2 r \,dr
		=	\\
		\displaystyle
		=	 
		3\int_0^{1/3}r^{2j+d-1} dr = \frac{3}{2j+d}\left(\frac{1}{3}\right)^{2j+d}< 2^{-2j}.
		\end{array} 		
	\end{equation}
		For $j < N$ we use fact that $||\psi_0(E)||_{L^2(B(0,1))}$ is continuous function on compact $W$ and, since $N$ depends only on $W$,  
	we get that there is a constant $\rho_1 = \rho_1(W,d)$ such that
	\begin{equation}
		 ||\psi_0(E)||_{L^2(B(0,1/3))} \leq \rho_1 2^{-j}.
	\end{equation}
		Since $v$ has support in $B(0,1/3)$ from (\ref{analytic_eq}) we get that
	\begin{equation}\label{l2sol}
		||\psi(E) - \psi_0(E)||_{L^2(B(0,1))}\leq \rho_1 2^{-j}||v||_{L^\infty(D)}||(-\Delta+v-E)^{-1}||_{L^2(D)}.
	\end{equation}
		Note that $\psi(E) - \psi_0(E)$ is the solution of equation (\ref{eq}) in 
	$D' = B(0,1)\setminus B(0,1/3)$ with potential $v_0 \equiv 0$ and boundary condition $\psi|_{r=1} = 0.$
		From Lemma \ref{Lemma_ring} we have that
	\begin{equation}\label{ck}
		\psi(E) - \psi_0(E) = \sum\limits_{0\leq i, 1\leq q\leq p_i} c_{iq}(E)\psi_{iq}(E) \  \text{ in $D'$} 
	\end{equation}
	for some $c_{iq}$, where
	\begin{equation}
		\psi_{iq}(E)(r,\omega) = R_{i}(\sqrt{E},r)f_{iq}(\omega).
	\end{equation}  
	Since $R_i(\sqrt{E},1) = 0$ 
	\begin{equation}
		\left.\frac{\partial R_i(\sqrt{E},r)}{\partial r}\right|_{r=1} = 
				\left.\frac{\partial \left({r^{\frac{d-2}{2}}}R_i(\sqrt{E},r)\right)}{\partial r}\right|_{r=1}.
	\end{equation}
	For $i \geq N$ from Lemma \ref{Bessel} we have that	
		\begin{equation}
			\begin{array}{c}
			\displaystyle
				\left|
				\frac{ \left.\frac{\partial R_i(\sqrt{E},r)}{\partial r}\right|_{r=1} }{Y_{\alpha}(\sqrt E)J_{\alpha}(\sqrt E)}  
				\right|
				=
				|E|^{1/2}
				\left|
					\frac{Y'_{\alpha}(\sqrt E)}{Y_{\alpha}(\sqrt E)}
					-
					\frac{J'_{\alpha}(\sqrt E)}{J_{\alpha}(\sqrt E)}
				\right| \leq
			\\\displaystyle
				\leq
				6 |E|^{1/2}
				\left(
						\frac
						{(|E|^{1/2}/2)^{-\alpha-1}\Gamma(\alpha+1)}   
						{(|E|^{1/2}/2)^{-\alpha}\Gamma(\alpha)}  
				+
						\frac
						{(|E|^{1/2}/2)^{\alpha-1}\Gamma(\alpha+1)}   
						{(|E|^{1/2}/2)^{\alpha}\Gamma(\alpha)} 
				\right) 
				=
				6\alpha,
			\end{array}
		\end{equation}
		\begin{equation}\label{ylemma1}
			\begin{aligned}
			\left(
			\frac
			{||r^{-\frac{d-2}{2}}Y_\alpha(\sqrt{E}r)||_{L^2(\{1/3 <|x|<2/5\})}}
			{|Y_{\alpha}(\sqrt E)|} 
			\right)^2
			&\geq
			\int_{1/3}^{2/5} 
				\left(
						\frac{1}{3}
						\frac
						{(|E|^{1/2}r/2)^{-\alpha}\Gamma(\alpha)}   
						{(|E|^{1/2}/2)^{-\alpha}\Gamma(\alpha)}  
				\right)^2 r\, dr \\
				&\geq 
				\left(\frac{2}{5} - \frac{1}{3}\right)\frac{1}{3} \left(\frac{1}{3} (5/2)^\alpha \right)^{2},
			\end{aligned}
		\end{equation}
		\begin{equation}\label{jlemma1}
			\begin{aligned}
			\left(
			\frac
			{||r^{-\frac{d-2}{2}}J_\alpha(\sqrt{E}r)||_{L^2(\{1/3 <|x|<2/5\})}}
			{|J_{\alpha}(\sqrt E)|} 
			\right)^2
			&\leq
			\int_{1/3}^{2/5} 
				\left(
						3
						\frac
						{(|E|^{1/2}r/2)^{\alpha}\Gamma(\alpha)}   
						{(|E|^{1/2}/2)^{\alpha}\Gamma(\alpha)}  
				\right)^2 r\, dr \\
				&\leq 
				\left(\frac{2}{5} - \frac{1}{3}\right)\frac{1}{3} \left(3 (2/5)^\alpha\right)^{2},
			\end{aligned}
		\end{equation}
		where $\alpha = i + \frac{d-2}{2}.$		Since $N > 3$ we have that $\alpha > 3$. 
		Using (\ref{ylemma1}) and (\ref{jlemma1}) we get that
		
\begin{equation}
		\frac
		{||\psi_{iq}(E)||_{L^2(\{1/3 <|x|<2/5\})}}
		{\left|Y_{\alpha}(\sqrt E)J_{\alpha}(\sqrt E)\right|}
		\geq
		\left(\Big(\frac{2}{5} - \frac{1}{3}\Big)\frac{1}{3}\right)^{1/2}
		\left(\frac{1}{3}(5/2)^\alpha - 3 (2/5)^\alpha\right)  
		\geq \frac{1}{1000}(5/2)^{\alpha}.
\end{equation}
	For $i \geq N$ we get that
	\begin{equation}
		\left|\left.\frac{\partial R_i(\sqrt{E},r)}{\partial r}\right|_{r=1}\right| \leq 1000\alpha(5/2)^{-\alpha}||\psi_{iq}(E)||_{L^2(\{1/3 <|x|<1\})}.
	\end{equation}		
	For $i < N$ we use the fact that 
	$\left|\left.\frac{\partial R_i(\sqrt{E},r)}{\partial r}\right|_{r=1}\right|/||\psi_{iq}(E)||_{L^2(\{1/3 <|x|<1\})}$ 
	is continuous function on compact $W$ and 
	get that for any $i\geq 0$ there is a constant $\rho_2 = \rho_2(W,d)$ such that	
	\begin{equation}\label{Lambda1}
		\left|\left.\frac{\partial R_i(\sqrt{E},r)}{\partial r}\right|_{r=1}\right| \leq \rho_2 \, 2^{-i}||\psi_{iq}(E)||_{L^2(\{1/3 <|x|<1\})}.
	\end{equation}
	Proceeding from (\ref{ck}) and using the Cauchy–Schwarz inequality we get that
	\begin{equation}\label{Lambda2}
		|c_{iq}(E)| = 
		\left|
		\frac
		{\Big< \psi(E) - \psi_0(E), \psi_{iq}(E) \Big>_{L^2(\{1/3 <|x|<1\})}}
		{||\psi_{iq}(E)||^2_{L^2(\{1/3 <|x|<1\})}} 
		\right|
		\leq
		\frac
		{||\psi(E) - \psi_0(E)||_{L^2(B(0,1))}}
		{||\psi_{iq}(E)||_{L^2(\{1/3 <|x|<1\})}}.
	\end{equation}
	Taking into account
	\begin{equation}
		\left\langle \Lambda_{v,E} f_{jp},f_{iq}\right\rangle
		=
		\left\langle  
		\left.\frac
		{\partial(\psi(E)-\psi_0(E))}
		{\partial\nu}\right|_{\partial D}, 
		f_{iq}
		\right\rangle
		=
		 c_{iq}(E)\left.\frac{\partial R_i(\sqrt{E},r)}{\partial r}\right|_{r=1}
	\end{equation}
	and combining (\ref{Lambda1}) and (\ref{Lambda2}) we obtain that
		\begin{equation}\label{l2sol2}
			\left|\left\langle \Lambda_{v,E} f_{jp},f_{iq}\right\rangle\right| \leq
			\rho_2 2^{-i}||\psi(E) - \psi_0(E)||_{L^2(B(0,1))}.
		\end{equation}
	From (\ref{l2sol})	and (\ref{l2sol2}) we get (\ref{eqlemma1}).

	For the general case we consider two compacts 
	\begin{equation}
		W_\pm = W \cap \left\{z \ | \ \pm\mbox{Im}z \geq 0 \right\}.
	\end{equation}
	Note that $\frac{J_{j+\frac{d-2}{2}}(\sqrt{E}r)}{J_{j+\frac{d-2}{2}}(\sqrt{E})}$ and
	 $\frac{Y_{j+\frac{d-2}{2}}(\sqrt{E}r)}{Y_{j+\frac{d-2}{2}}(\sqrt{E})}$
	have removable singularity in $E=0$ or, more precisely,
	\begin{equation}
		\begin{aligned}
		\frac{J_{j+\frac{d-2}{2}}(\sqrt{E}r)}{J_{j+\frac{d-2}{2}}(\sqrt{E})} &\longrightarrow r^{j+\frac{d-2}{2}}, \\\displaystyle 
		\frac{Y_{j+\frac{d-2}{2}}(\sqrt{E}r)}{Y_{j+\frac{d-2}{2}}(\sqrt{E})} &\longrightarrow r^{-j-\frac{d-2}{2}}
		\\ 
		\text{ as } E &\longrightarrow 0.
		\end{aligned}
	\end{equation}
 	Considering the limit as $E \rightarrow 0$ we get that (\ref{l2sol}), (\ref{l2sol2}) and consequently (\ref{eqlemma1}) are valid for $W_\pm$. 
 	To complete proof we can take $\rho = \max\{\rho_+, \rho_-\}$. 
\end{Proof}
 	
\noindent
	{\bf Remark 3.2.}
 	From  (\ref{psi_0}) and (\ref{analytic_eq}) we get that
 	\begin{equation}
 		\left\langle \Lambda_{v,E} f_{jp},f_{iq}\right\rangle\ \text{ is holomorphic function in } W.
 	\end{equation}


\section{A fat metric space and a thin metric space}

\begin{Def}
	Let $(X,dist)$ be a metric space and $\epsilon > 0$. We say that a set $Y \subset X$ is an $\epsilon$-net
	for $X_1 \subset X$ if for any $x\in X_1$ there is $y \in Y$ such that $dist(x,y)\leq\epsilon.$ We call 
	$\epsilon$-entropy of the set $X_1$ the number $\mathcal H_\epsilon (X_1) := \log_2
	\min\{|Y|: Y$  is an $\epsilon$-net fot $X_1 \}.$ 
	
	A set $Z \subset X$ is called $\epsilon$-discrete if for any distinct $z_1, z_2 \in Z$, we have 		  
	$dist(z_1,z_2)\geq\epsilon$. We call $\epsilon$-capacity of the set $X_1$ the number $\mathcal C_\epsilon :=
	 \log_2	\max\{|Z|: Z \subset X_1 $ and $Z$ is $\epsilon$-discrete$\}.$
\end{Def}

	The use of $\epsilon$-entropy and $\epsilon$-capacity  to derive properties of mappings between metric
spaces goes back to Vitushkin and Kolmogorov (see \cite{Kolmogorov1959} and references therein). One notable
application was Hilbert’s 13th problem (about representing a function of several variables as
a composition of functions of a smaller number of variables). In essence, Lemma \ref{Lemma2} and Lemma \ref{Lemma3} are parts of the Theorem XIV and the Theorem XVII in \cite{Kolmogorov1959}. 

\begin{Lemma}\label{Lemma2}
	Let $d \geq 2$ и $m > 0$. For $\epsilon, \beta > 0$, consider the real metric space
	$$
		X_{m\epsilon\beta} = \{f \in C^{m}(D)\ |\   \mbox{supp}\, f \subset B(0,1/3),\  
		||f||_{L^\infty(D)} \leq \epsilon,\ 	||f||_{C^m(D)}\leq \beta\},
	$$
	with the metric induced by $L^{\infty}$. Then there is a $\mu>0$ such that for any $\beta > 0$ and 			
	$\epsilon \in (0, \mu\beta)$, there is an
	$\epsilon$-discrete set  $Z \subset X_{m\epsilon\beta}$ with at least $\exp\Big(2^{-d-1}(\mu\beta/\epsilon)^{d/m}\Big)$ elements.
\end{Lemma}

Lemma \ref{Lemma2} was also formulated and proved in \cite{Mandache2001}.
	
\begin{Lemma}\label{Lemma3}
	For the interval $I = [a,b]$ with $a < b$ and $\gamma>0$ consider ellipse $W_{I,\gamma} \in \mathbb{C}$
	\begin{equation}\label{rectangle}
		W_{I,\gamma} = \{ \frac{a+b}{2}+ \frac{a-b}{2}\cos z \ |\ |Im\, z| \leq \gamma \}.
	\end{equation}  	
	Then there is a constant $\nu = \nu(C,\gamma)>0$, such that for every $\delta\in(0,e^{-1})$, there is
	a $\delta$-net for the space functions on $I$ with $L^{\infty}$-norm, having holomorphic continuation to $W_{I,\gamma}$ 
	 with module bounded above on $W_{I,\gamma}$ by the constant $C$, 
	with at most $\exp(\nu (\ln\delta^{-1})^2)$ elements.    
\end{Lemma}

\begin{Proof} {\it Lemma \ref{Lemma3}. }
	Theorem XVII in \cite{Kolmogorov1959} provides asymptotic behaviour of the entropy of this space with respect to $\delta\rightarrow 0$. 
	Here we get upper estimate of it. 
	Suppose $g(z)$ is holomorphic function in $W_{I,\gamma}$ with module bounded above by the constant $C$.
	Consider the function $f(z) = g(\frac{a+b}{2}+ \frac{a-b}{2}\cos z)$. By the choise of $W_{I,\gamma}$ we get that 
	$f(z)$ is $2\pi$-periodic holomorphic function in the stripe $|\mbox{Im}\,z|\leq \gamma$. Then for any integer $n$
	\begin{equation}
		|c_n| = \left|\int_0^{2\pi} e^{inx}f(x) dx\right| \leq \int_0^{2\pi} e^{-|n|\gamma}C dx \leq 2\pi C e^{-|n|\gamma}.
	\end{equation}
	Let $n_\delta$ be the smallest natural number such that $2\pi C e^{-n\gamma} \leq 6 \pi^{-2} (n+1)^{-2}\delta$ for any $n \geq n_\delta$.
	Taking natural logarithm and using $\ln\delta^{-1} \geq 1$, we get that
	\begin{equation}\label{ndelta}
		n_\delta \leq C'\ln\delta^{-1}, 
	\end{equation}
	where $C'$ depends only on $C$ and $\gamma$. 
	We denote $\delta' = {3} \pi^{-2} (n_\delta+1)^{-2}\delta$.
	Consider the set 
	\begin{equation}
		Y_\delta = \delta'\mathbb Z \bigcap [-2\pi C, 2\pi C] + i\cdot\delta'\mathbb Z \bigcap [-2\pi C, 2\pi C].
	\end{equation} 
	Using (\ref{ndelta}), we have that
	\begin{equation}\label{nydelta}
		|Y_\delta| = \left(1 + 2[2\pi C/\delta']\right)^2 \leq C''\delta^{-2}\ln^4\delta^{-1},
	\end{equation} 
	with $C''$ depending only on $C$ and $\gamma$. 
	We set 
	\begin{equation}
		Y = 
		\left\{ 
		\sum\limits_{n = 0}^{\infty} d_n \cos\left(n\arccos{\frac{x-\frac{a+b}{2}}{\frac{a-b}{2}}}\right) \ \left|\ 
		d_n \in Y_\delta \text{ for $n\leq n_\delta$, $d_n = 0$ otherwise} \right. 
		\right\}.
	\end{equation}
	For given $f(z)$ in case of  $n\leq n_\delta$ we take $d_n$ to be one of the closest elements of $Y_\delta$ to $c_n$. 
	Since $|c_n|\leq 2\pi C$, this ensures $|c_n-d_n| \leq 2\delta'$. For  $n > n_\delta$ we take $d_n = 0.$ We have then
	\begin{equation}
		|c_n-d_n| \leq 6 \pi^{-2} (n+1)^{-2}\delta.
	\end{equation}
	For $n>n_\delta$ this is true by the construction of $n_\delta$, otherwise by the choise of $\delta'$. 
	Since $f(x)$ is $2\pi$-periodic even function, we get $g_Y(x) \in Y$ such that 
	\begin{equation}
		||g(x)-g_Y(x)||_{L^{\infty}(a,b)} \leq \sum\limits_{n=0}^{\infty} |c_n-d_n|
		\leq 6 \pi^{-2}\delta \sum\limits_{n=1}^{\infty} \frac{1}{n^2} = \delta.
	\end{equation} 
		We have that $|Y| = |Y_\delta|^{n_\delta}$. Taking into account (\ref{ndelta}),(\ref{nydelta}) and $\ln\delta^{-1} \geq 1$, we get  
	\begin{equation}
		|Y|\leq (C''\delta^{-2}\ln^4\delta^{-1})^{C'\ln\delta^{-1}} \leq 
		\exp\left( C'''\ln\delta^{-1}C'\ln\delta^{-1} \right)\leq \exp(\nu (\ln\delta^{-1})^2).
	\end{equation} 
\end{Proof}

\noindent
{\bf Remark 4.1.} The assertion is valid even in the case of $a = b$. As $\delta$-net we can take
	\begin{equation}
		Y =\frac{\delta}{2}\mathbb{Z}\bigcap[-C,C] + i\cdot \frac{\delta}{2}\mathbb{Z}\bigcap[-C,C].
	\end{equation}

	Consider an operator $A: H^{-s}(S^{d-1})\rightarrow H^{s}(S^{d-1})$. We denote its matrix elements in the 
	basis $\{f_{jp}\}$ by $a_{jpiq} = \left\langle Af_{jp}, f_{iq}\right\rangle$. From \cite{Mandache2001} we have that
	\begin{equation}\label{AHS}
		||A||_{H^{-s}\rightarrow H^{s}}\leq 4 \sup\limits_{j,p,i,q}(1+\max(j,i))^{2s+d}|a_{jpiq}|.
	\end{equation}
	Consider system $S = \bigcup\limits_{j=1}\limits^{K} I_j$ of $\sigma$-regular intervals. 
 	We introduce the Banach space
		$$
			X_{S,s} = \left\{ \Big(a_{jpiq}(E)\Big)\ | \ 
			 \left\|\Big(a_{jpiq}(E)\Big)\right\|_{X_{S,s}} := \sup_{j,p,i,q} \left((1 + \max(j,i))^{2s+d}\sup_{E\in S}|a_{jpiq}(E)|\right) < \infty \right\}.
		$$
 	Denote by $B^{\infty}$ the ball of centre $0$ and radius $2\sigma/3$ in $L^\infty(B(0,1/3))$. We identify in the sequel 
 	an operator $A(E): H^{-s}(S^{d-1})\rightarrow H^{s}(S^{d-1})$ with its matrix $\Big(a_{jpiq}(E)\Big)$. Note that the estimate (\ref{AHS}) 
 	implies that 
 	\begin{equation}\label{xssajpiq}
 		\sup_{E \in S}\left\|A(E)\right\|_{H^{-s}\rightarrow H^{s}} \leq 4 \left\|\Big(a_{jpiq}(E)\Big)\right\|_{X_{S,s}}.
 	\end{equation} 

We consider operator $\Lambda_{v,E}$ from Lemma \ref{Lemma1} as
\begin{equation}
	\Lambda: B^{\infty} \rightarrow \left\{ \Big(a_{jpiq}(E)\Big) \right\},
\end{equation}
where $a_{jpiq}(E)$ are matrix elements in the basis $\{f_{jp}\}$ of operator $\Lambda_{v,E}$.

\begin{Lemma}\label{Lemma4}
$\Lambda$ maps $B^{\infty}$ into $X_{S,s}$ for any $s$. There is a constant $0 < \eta = \eta(S,s,d)$, such that 
for every $\delta \in (0,e^{-1})$, there is a $\delta$-net $Y$ for $\Lambda(B^{\infty})$ in $X_{S,s}$  
with at most $\exp(\eta(\ln\delta^{-1})^{2d})$ elements.
\end{Lemma}

\begin{Proof} {\it Lemma \ref{Lemma4}. }
	For simplicity we give first a proof in case of $S$ consists of only one $\sigma$-regular interval $I$. 
	From (\ref{rectangle}) we take $W_I = W_{I,\gamma}$, where constant $\gamma>0$ is such as for any $E\in W_I$ 
	there is $E_I$ in $I$  such as $|E-E_I|<\sigma/6$. From (\ref{lambdaE}) we get that 
	\begin{equation}
		|E-\lambda| \geq |E_I-\lambda| - |E-E_I| \geq 5\sigma/6,
	\end{equation} 		 
	with $\lambda$ being Dirichlet eigenvalue for operator $-\Delta$ in $D$ which is closest to $E$.
	Then for potential $v \in B^{\infty}$ and $E\in W_I$ we have that
	\begin{equation}	
	 ||(-\Delta+v-E)^{-1}||_{L^2(D)}\leq (|\lambda-E| - 2\sigma/3)^{-1} \leq (5\sigma/6 - 2\sigma/3)^{-1} = 6/\sigma   
	\end{equation}
	 and
	 \begin{equation}
	 	||v||_{L^{\infty}(D)}||(-\Delta+v-E)^{-1}||_{L^2(D)} \leq (2\sigma/3)(6/\sigma) = 4, 
	 \end{equation} 
	where $(-\Delta+v-E)^{-1}$ is considered with the Dirichlet boundary condition.
	 We obtain from Lemma \ref{Lemma1} that  
	\begin{equation}\label{ajpkqbound}
		|a_{jpiq}(E)| \leq  4\rho\,2^{-\max(j,i)},
	\end{equation}
	where $\rho = \rho(W_I,d)$.
	Hence $||(a_{jpiq}(E))||_{X_{S,s}}\leq \sup_l(1+l)^{2s+d}4\rho\,2^{-l}<\infty$ for any $s$ and $d$ and so the first assertion of the 
	Lemma \ref{Lemma4} is proved.

		Let $l_{\delta s}$ be the smallest natural number such that $(1+l)^{2s+d}4\rho\,2^{-l} \leq \delta$ for any $l \geq l_{\delta s}.$
		Taking natural logarithm and using $\ln\delta^{-1} \geq 1$, we get that
		\begin{equation}\label{lsmax}
		l_{\delta s} \leq C'\ln\delta^{-1},
		\end{equation}
		where the constant $C'$ depends only on $s$, $d$ and $I$.
		 Denote $Y_{jpiq}$ is $\delta_{jpiq}$-net	from Lemma \ref{Lemma3} with constant $C = \sup_l(1+l)^{2s+d}4\rho\,2^{-l},$ 
		where $\delta_{jpiq} = (1+\max(j,i))^{-2s-d}\delta$. We set
		\begin{equation}
			Y = \left\{
			(a_{jpiq}(E)) \ |\ a_{jpiq}(E)\in Y_{jpiq} 
			\text{ for } \max(j,i)\leq l_{\delta s},
			\ a_{jpiq}(E) = 0 \text{ otherwise} 
			 \right\}.
		\end{equation}
		For any $(a_{jpiq}(E)) \in \Lambda(B^{\infty})$ there is an element $(b_{jpiq}(E)) \in Y$ such that
	
		\begin{equation}
			(1+\max(j,i))^{2s+d}|a_{jpiq}(E)-b_{jpiq}(E)| \leq (1+\max(j,i))^{2s+d}\delta_{jpiq} = \delta,
		\end{equation}
		in case of $\max(j,i) \leq l_{\delta s}$ and
			\begin{equation}
				(1+\max(j,i))^{2s+d}|a_{jpiq}(E)-b_{jpiq}(E)|		
				\leq	(1+\max(j,i))^{2s+d}2\rho\,2^{-\max(j,i)} \leq \delta,
			\end{equation}
		otherwise.
		
			It remains to count the elements of $Y$. Using again the fact that $\ln\delta^{-1} \geq 1$ and (\ref{lsmax})
			we get for $\max(j,i)\leq l_{\delta s}$  
			\begin{equation}
				|Y_{jpiq}| \leq \exp(\nu (\ln\delta_{jpiq}^{-1})^2) \leq \exp(\nu' (\ln\delta^{-1})^2).
			\end{equation}
			From \cite{Mandache2001} we have that $n_{\delta s} \leq 8(1+l_{\delta s})^{2d-2},$ where 
			$n_{\delta s}$ is the number of four-tuples $(j,p,i,q)$ with $\max(j,i) \leq l_{\delta s}$. 
			Taking $\eta$ to be big enough we get that
			\begin{equation}
				\begin{aligned}
				|Y| &\leq \left(\exp(\nu' (\ln\delta^{-1})^2)\right)^{n_{\delta s}} \\
						&\leq \exp\left(\nu' (\ln\delta^{-1})^2 {8(1+C'\ln\delta^{-1})^{2d-2}}\right)\\
						&\leq \exp \left(\eta(\ln\delta^{-1})^{2d}\right).
				\end{aligned}	
			\end{equation}
			For $S = \bigcup\limits_{j=1}\limits^{K} I_j$ assertion follows immediately, taking $\eta$ to be in $K$ times more and 
			$Y$ as composition  $(Y_1,\ldots,Y_K)$ of $\delta$-nets for each interval.  
\end{Proof}

\section{Proofs of the main results}
In this section we give proofs of Theorem \ref{main} and Theorem \ref{complex}.\\

\begin{Proof} {\it Theorem \ref{main}.}
	Take $v_0 \in L^{\infty}(B(0,1/3))$, $||v_0||_{L^\infty(D)}\leq \sigma/3$ and $\epsilon \in (0,\sigma/3).$
	By Lemma \ref{Lemma2}, the set $v_0 + X_{m\epsilon\beta}$ has an $\epsilon$-discrete subset $v_0 + Z$. 
	Since for  $\epsilon \in (0,\sigma/3)$ we have $v_0+X_{m\epsilon\beta} \subset B^{\infty},$ where $B^\infty$ is 
	the ball of centre $0$ and radius $2\sigma/3$ in $L^\infty(B(0,1/3))$.
	The set $Y$
	constructed in Lemma \ref{Lemma4} is also $\delta$-net for $\Lambda(v_0+X_{m\epsilon\beta})$. 
	We take $\delta$ such that $8\delta = \exp\left(-\epsilon^{-\frac{1}{2m}} \right)$. Note that inequalities of (\ref{maineq}) follow from
	\begin{equation}\label{ZplusY}
		|v_0 + Z|>|Y|.
	\end{equation}
	In fact, if $|v_0 + Z|>|Y|$, then there are two potentials $v_1, v_2 \in v_0 + Z$ with images under $\Lambda$
	in the same $X_{S,s}$-ball radius $\delta$ centered at a point of $Y$, so we get from (\ref{xssajpiq})
	\begin{equation}
		\sup\limits_{E\in S} ||\Phi_1(E)-\Phi_2(E)||_{H^{-s} \rightarrow H^{s}} 
		\leq 
			4||\Lambda_{v_1,E}-\Lambda_{v_2,E}||_{X_{S,s}} 
		\leq 
		8\delta = \exp\left(-\epsilon^{-\frac{1}{2m}} \right).
	\end{equation}
	It remains to find $\beta$ such as (\ref{ZplusY}) is fullfiled.	By Lemma \ref{Lemma4}
	\begin{equation} \label{eqy}
		|Y|\leq \exp\left(\eta\left( \ln8 +\epsilon^{-\frac{1}{2m}} \right)^{2d} \right)
		\leq \max\Big(  \exp\left( (2\ln8)^{2d}\eta\right),  \exp\left( 2^{2d}\eta\epsilon^{-d/m}\right) \Big).
	\end{equation}
	Now we take 
	\begin{equation}\label{eqbeta}
		\beta > \mu^{-1} \max\left(\sigma/3, \eta^{m/d}2^{3m}, \frac{\sigma}{3}\eta^{m/d}2^{m}(2\ln8)^{2m} \right) 
	\end{equation}
	This fulfils requirement $\epsilon<\mu\beta$ in Lemma \ref{Lemma2}, which gives
	\begin{equation}
		\begin{array}{c}\displaystyle
			|v_0+Z| = |Z| \geq \exp\Big(2^{-d-1}(\mu\beta/\epsilon)^{d/m}\Big)\stackrel{(\ref{eqbeta})}{>}\\\displaystyle
			>\max\Big(\exp\left(2^{-d-1}(\eta^{m/d}2^{3m}/\epsilon)^{d/m}\right),  
									\exp\left(2^{-d-1}(\eta^{m/d}2^{m}(2\ln8)^{2m})^{d/m}\right)\Big) \stackrel{(\ref{eqy})}{\geq} |Y|.
		\end{array}
	\end{equation}
\end{Proof}

\begin{Proof} {\it Theorem \ref{complex}.} 
	In a similar way with the proof of Theorem 2 of \cite{Mandache2001} we obtain that
	\begin{equation}\label{jizero}
		\left\langle \left(\Phi_{mn}(E) - \Phi_0(E)\right)f_{jp}, f_{iq}\right\rangle = 0
	\end{equation}  
	for $j,i \leq \left[\frac{n-1}{2}\right]$. The only difference is that instead of the operator $-\Delta$ we consider
	the operator $-\Delta - E$.
	From (\ref{AHS}), (\ref{ajpkqbound}) and (\ref{jizero}) we get
	\begin{equation}
	 ||\Phi_{mn}(E) - \Phi_0(E)||_{H^{-s} \rightarrow H^{s}} \leq 16\rho\sup_{l\geq n/2}(1+l)^{2s+d}2^{-l} \leq c'2^{-n/4}.
	\end{equation}
	The fact that $||v_{mn}||_{C^m(D)}$ is bounded as $n \rightarrow \infty$ is also a part of Theorem 2 of \cite{Mandache2001}.
\end{Proof}

\section{Bessel functions}

In this section we prove Lemma \ref{Lemma_solution}, Lemma \ref {Lemma_ring} and Lemma \ref{Bessel} about the Bessel functions. 
 Consider the problem of finding solutions of the form $\psi(r,\omega) = R(r)f_{jp}(\omega)$
of equation (\ref{eq}) with $v\equiv 0$ . We have that
\begin{equation}
	\Delta = \frac{\partial^2}{(\partial r)^2} + (d-1)r^{-1} \frac{\partial}{\partial r} + r^{-2}\Delta_{S^{d-1}},   
\end{equation}
where $\Delta_{S^{d-1}}$ is  Laplace-Beltrami operator on $S^{d-1}$. We have that
\begin{equation}
	\Delta_{S^{d-1}}f_{jp} = -j(j+d-2)f_{jp}.
\end{equation}
Then we have the following equation for $R(r)$:
\begin{equation}
	-R'' - \frac{d-1}{r}R' + \frac{j(j+d-2)}{r^2}R = ER.
\end{equation}
Taking $R(r) = r^{-\frac{d-2}{2}}\tilde{R}(r)$, we get
\begin{equation}
	r^2\tilde{R}'' + r\tilde{R}' + \left(Er^2 - \left(j+\frac{d-2}{2}\right)^2\right)\tilde{R} = 0.
\end{equation}
This equation is known as Bessel's equation. For $E =k^2\neq 0$ it has two linearly independent solutions 
$J_{j+\frac{d-2}{2}}(kr)$ and $Y_{j+\frac{d-2}{2}}(kr),$ where 
\begin{equation}\label{jalphaz}
	J_\alpha(z) = \sum\limits_{m=0}\limits^{\infty}
		\frac{(-1)^m(z/2)^{2m+\alpha}}{\Gamma(m+1)\Gamma(m+\alpha + 1)},
\end{equation} 
\begin{equation}
	Y_\alpha(z) = \frac{J_\alpha(z)\cos\pi\alpha  - J_{-\alpha}(z)}{\sin\pi\alpha} \text{ for  $\alpha \notin \mathbb{Z}$,}
\end{equation}
and 
\begin{equation}
	Y_\alpha(z) = \lim\limits_{\alpha'\rightarrow \alpha}Y_{\alpha'}(z) \text{ for  $\alpha \in \mathbb{Z}$.}
\end{equation}
The following Lemma is called the Nielsen inequality. A proof can be found in \cite{Watson1944} 
\begin{Lemma}\label{Nielsen}
	\begin{equation}
	\begin{array}{l}\displaystyle
		J_\alpha(z) = \frac{\left(z/2\right)^\alpha}{\Gamma(\alpha+1)}(1+\theta),
		\\\displaystyle
		|\theta| < \exp\left(\frac{|z|^2/{4}}{|\alpha_0+1|}\right) - 1,
		\end{array}
	\end{equation}	
		where $|\alpha_0+1|$ is the least of numbers $|\alpha+1|,|\alpha+2|,|\alpha+3|,\ldots$ .
\end{Lemma}
Lemma \ref{Nielsen} implies that $r^{-\frac{d-2}{2}}J_{j+\frac{d-2}{2}}(kr)$ has removable singularity at $r = 0$. 
Using the boundary conditions $R(1) = 1$ and $R(1) = 0$, we obtain assertions of Lemma \ref{Lemma_solution} and Lemma \ref{Lemma_ring}, respectively.

\begin{Proof} {\it Lemma \ref{Bessel}}
	Formula (\ref{Jn}) follows immediately from Lemma \ref{Nielsen}.  We have from \cite{Watson1944} that
	\begin{equation}
		J_\alpha'(z) = J_{\alpha-1}(z) - \frac{\alpha}{z} J_\alpha(z).
	\end{equation}
	Further, taking $\alpha$ big enough we get
	\begin{equation}
		|J_\alpha'(z)| \leq |J_{\alpha-1}(z)| + |\frac{\alpha}{z} J_\alpha(z)| 
		\leq \frac{3}{2} \frac{\left(|z|/2\right)^{\alpha-1}}{\Gamma(\alpha)} + \frac{3\alpha}{2|z|} \frac{\left(|z|/2\right)^{\alpha}}{\Gamma(\alpha+1)}
		\leq 3 \frac{\left(|z|/2\right)^{\alpha-1}}{\Gamma(\alpha)}.
	\end{equation}
	For $\alpha = n + 1/2$ we have $Y_\alpha = (-1)^{n+1}J_{-\alpha}$. Consider its series expansion, see (\ref{jalphaz}). 
	\begin{equation}
		J_{-\alpha}(z) = \sum\limits_{m=0}\limits^{\infty}
		\frac{(-1)^m(z/2)^{2m-\alpha}}{m!\,\Gamma(m-\alpha + 1)} = \sum\limits_{m=0}\limits^{\infty}c_m(z/2)^{2m-\alpha}.
	\end{equation}
	Note that $|c_{m}/c_{m+1}| = (m+1)|m-\alpha+1| \geq n/2$. As corollary we obtain that
	\begin{equation}\label{proofhn}
		\begin{array}{c}\displaystyle
		|Y_\alpha(z)| = \frac{(|z|/2)^{-\alpha}}{|\Gamma(-\alpha + 1)|}(1+\theta) = \frac{1}{\pi}(|z|/2)^{-\alpha}\Gamma(\alpha)(1+\theta),
		\\\displaystyle
			|\theta| \leq \sum\limits_{m=1}\limits^{\infty}\left(\frac{|z|^2}{2n}\right)^{2m}
			\leq \frac{|z|^2/2n}{1-|z|^2/2n}.
		\end{array}
	\end{equation} 
	For $\alpha = n$ we have from \cite{Watson1944} that
	\begin{equation}\label{bmcm}
	\begin{array}{c}\displaystyle
		Y_n(z) = \frac{2}{\pi}J_n(z)\ln\left(\frac{z}{2}\right) -
		\frac{1}{\pi} \sum\limits_{m=0}\limits^{n-1}\left(\frac{z}{2}\right)^{2m-n}\frac{(n-m-1)!}{m!} -\\
		\displaystyle-
		\frac{1}{\pi} \sum\limits_{m=0}\limits^{\infty}\frac{(-1)^m(z/2)^{2m+n}}{m!(m+n)!}
		\left(\frac{\Gamma'(m+1)}{\Gamma(m+1)}+ \frac{\Gamma'(m+n+1)}{\Gamma(m+n+1)}\right)=
		\\\displaystyle=
		\frac{2}{\pi}J_n(z)\ln\left(\frac{z}{2}\right) - \frac{1}{\pi} \sum\limits_{m=0}\limits^{n-1}\tilde{c}_m(z/2)^{2m-n}
		-\frac{1}{\pi} \sum\limits_{m=0}\limits^{\infty}b_m(z/2)^{2m+n}.
		\end{array}
	\end{equation}
	Using well-known equality $\Gamma'(x)/\Gamma(x) < \ln x$, $x > 1$, see \cite{Abramowitz1972}, 
	we get following estimation for the coefficients $b_m$ are defined in (\ref{bmcm}).  
	\begin{equation}\label{bmbm}
		|b_m| < \frac{\ln(m+1)+\ln(n+m+1)}{m!(n+m)!} < \frac{2(n+m)}{m!(n+m)!}< \frac{1}{m!}.
	\end{equation}
	Note also that $|\tilde{c}_{m}/\tilde{c}_{m+1}| = (m+1)(n-m-1) \geq n/2.$ Combining it with (\ref{bmcm}) and (\ref{bmbm}), we obtain that
	\begin{equation}\label{proofn}
		\begin{array}{c}\displaystyle
		|Y_n(z)| = \frac{1}{\pi}(|z|/2)^{-n}\Gamma(n)(1+\theta),
		\\\displaystyle
			|\theta| \leq 
			3\frac{(|z|/2)^{2n}|\ln(z/2)|}{\Gamma(n)}+
			\sum\limits_{m=1}\limits^{n-1}\left(\frac{|z|^2}{2n}\right)^{2m} 
			+\frac{(|z|/2)^{2n}}{\Gamma(n)}\sum\limits_{m=0}\limits^{\infty}\frac{(|z|/2)^{2m}}{m!}
			\leq 
			\\\displaystyle
			\leq 3\pi\frac{\max\left(1,(|z|/2)^{2n+1}\right)}{\Gamma(n)} + 
			\frac{|z|^2/2n}{1-|z|^2/2n} + \frac{(|z|/2)^{2n}e^{|z|^2/4}}{\Gamma(n)}.
		\end{array}
	\end{equation} 
	Formula (\ref{Yn}) follows from (\ref{proofhn}) and (\ref{proofn}). We have from \cite{Watson1944} that
	\begin{equation}
		Y_\alpha'(z) = Y_{\alpha-1}(z) - \frac{\alpha}{z} Y_\alpha(z).
	\end{equation}
	Taking $n$ big enough, we get that
	\begin{equation}
	\begin{array}{c}\displaystyle
		|Y_\alpha'(z)| \leq |Y_{\alpha-1}(z)| + |\frac{\alpha}{z} Y_\alpha(z)| 
		\leq 
		\\\displaystyle \leq
		\frac{3}{2\pi} \left( \left(|z|/2\right)^{-\alpha+1}{\Gamma(\alpha-1)} + \frac{\alpha}{|z|} {\left(|z|/2\right)^{\alpha}}{\Gamma(\alpha)}\right)
		\leq  \frac{3}{\pi}(|z|/2)^{-\alpha-1}\Gamma(\alpha+1).
		\end{array}
	\end{equation}
	Combining reqirements for $n$, stated above, we get that for any $n \geq N+1$ all inequalities of Lemma \ref{Bessel} are fullfiled, where $N$ such that
	\begin{equation}\left\{
		\begin{array}{l}\displaystyle
				N>3,
			\\\displaystyle
			\exp\left(\frac{C^2/{4}}{N+1}\right) - 1 \leq 1/2 ,
				\\\displaystyle
				3\pi\frac{\max\left(1,(C/2)^{2N+1}\right)}{\Gamma(N)} + 
			\frac{C^2}{2N-C^2} + \frac{(C/2)^{2N}e^{C^2/4}}{\Gamma(N)} \leq 1/2.
		\end{array}\right.
	\end{equation}
	\end{Proof}

\section*{Acknowledgments}
This work was fulfilled under the direction
of R.G.Novikov in the framework of an internship at Ecole Polytechnique.


\begin{thebibliography}{99}
\bibitem{Alessandrini1988}
 G.Alessandrini, 
 {\it Stable determination of conductivity by boundary measurements}, 
 Appl.Anal. 27 (1988) 153-172.


\bibitem{Gelfand1954}
I.M. Gelfand, {\it Some problems of functional analysis and algebra}, Proceedings of the
International Congress of Mathematicians, Amsterdam, 1954, pp.253-276.

\bibitem{Mandache2001}
 N. Mandache,
{ \it Exponential instability in an inverse problem for the Schr\"odinger equation}
 { Inverse Problems}. 17(2001) 1435–1444. 

\bibitem{Henkin1987}
G.M. Henkin and R.G. Novikov, 
{\it The $\bar{\partial}$-equation in the multidimensional inverse scattering problem}, 
Uspekhi Mat. Nauk 42(3) (1987), 93-152 (in Russian); English
Transl.: Russ. Math. Surv. 42(3) (1987), 109-180.

\bibitem{Watson1944}
G. N. Watson,  
{ \it A Treatise on the Theory of Bessel Functions.} 
Cambridge University Press, Cambridge, England; The Macmillan Company, New York, 1944.

\bibitem{Novikov1988}
R.G. Novikov, {\it Multidimensional inverse spectral problem for the equation $-\Delta \psi + (v(x) - Eu(x))\psi = 0$}
Funkt. Anal. Prilozhen. 22(1988) 11–22 (in Russian) (Engl. Transl.  Funct. Anal. Appl. 22(1988) 263–72).

\bibitem{Novikov2010}
R.G. Novikov,
{\it New global stability estimates for the Gel’fand-Calderon inverse problem,} e-print arXiv:1002.0153.

\bibitem{NS}
R. Novikov and M. Santacesaria,
{\it A global stability estimate for the Gel’fand- Calderon
inverse problem in two dimensions,} e-print arXiv: 1008.4888.

\bibitem{Bukhgeim2008}
 A. L. Bukhgeim, 
 {\it Recovering a potential from Cauchy data in the two-dimensional
case}, J. Inverse Ill-Posed Probl. 16, 2008, no. 1, 19–33.













\bibitem{Kolmogorov1959}
A.N. Kolmogorov,V.M. Tikhomirov 
{\it $\epsilon$-entropy and $\epsilon$-capacity in functional spaces} Usp. Mat. Nauk 14(1959)
3–86 (in Russian) (Engl. Transl.  Am. Math. Soc. Transl. 17 (1961) 277–364)

\bibitem{Abramowitz1972}
M. Abramowitz, I.A. Stegun,(Eds.). 
{\it Psi (Digamma) Function.} §6.3 in Handbook of Mathematical Functions with Formulas, Graphs, and Mathematical Tables, 9th printing. New York: Dover, pp. 258-259, 1972

\bibitem{Cristo2003}
M. Di Cristo and L. Rondi
{\it Examples of exponential instability for inverse
inclusion and scattering problems}
{ Inverse Problems}. 19 (2003) 685–701.

\end{thebibliography}
\end{document}